# The Sheaf of the Groups Formed by Topological Generalized Group over Topological spaces

Hatice ASLAN[1,*], Hakan EFE[2]

[1] *Firat University, Faculty of Science, Department of Mathematics,*
*23119, Elazig, Turkey, email: haslan@firat.edu.tr*

[2]*Gazi University, Faculty of Science, Depatment of Mathematics,*
*06500 Teknikokullar Ankara, Turkey, email: hakanefe@gazi.edu.tr*

**ABSTRACT**

In the present paper, we show how to construct an algebraic sheaf by means of the topological generalized group defined by Molaei in [16] by considering both homotopy and sheaf theory.

**Keywords:** Generalized group, Topological generalized group, Whitney sum, Sheaf.
**2010 Mathematics Subject Classification:** 22A05; 14H30; 46M20.

## 1. INTRODUCTION AND PRELIMINARIES

Generalized groups were deduced from a geometrical defined by Molaei [16] in 1999. It is an interesting and fascinating extension of groups. Molaei's generalized groups established the uniqueness of the identity element of each element in a generalized group and where the identity element is not unique for each element. With this property, every group is a generalized group. This new concept studied in terms of algebraic, topological and differentiable in large various areas of mathematics [2, 13–18].

Generalized group is an algebraic structure which has a deep physical background in the unified guage theory. The unified theory has a direct relation with the geometry of space. It describes particles and their interactions in a quantum mechanical manner and the geometry of the space-time through which they are moving. Currently, the most promising is super-string theory in which the so called elementary particles are described as vibration of tiny (planck-length) closed loops of strings. In this theory the classical law of physics, such as electromagnetism and general relativity, are modified at time distances comparable to the length of the string. This notion of 'quantum space-time' is the goal of unified theory of physical forces.

Therefore the unified theory offers a new insight into the structure, order and measures of the quantum world of the entire universe. It is known that unified theories are based on the geometry of a space and the metric can determine the geometry [8]. Because of this physical forces mathematicians and physicists have been trying to construct a suitable unified theory for unified gauge theory, twistor theory, isotopies theory and so on. Now it is known that generalized groups are tools for constructions in unified geometric theory and electroweak theories which is essentially structured on Minkowskian axioms and gravitational theories are constructed on Riemannian axioms.

Furthermore this kind of structure appears in genetic codes. Generalized groups have been applied to DNA analysis by transforming the set of DNA sequences to generalized group in [1].

Another important concept in this present paper are sheaves which were originally introduced by Leray [12] in 1946. The modified definition of sheaves now used was given by Lazard, and appeared first in the Cartan Sem. [5] 1950-51. Sheaf theory provides a language for the discussion of geometric objects of many different kinds. At present it finds its main applications in topology and (more especially) in modern algebraic geometry, where it has been used with great success as a tool in the solution of several longstanding problems.

C. Yıldız constructed an algebraic sheaf by means of the topological group in [23]. This is our motivation for costructing a sheaf by the means of the topological generalized group in this paper. We replace topological group with topological generalized group construct an algebraic sheaf by means of the topological generalized group defined by Molaei in [16].

This section of the paper is devoted to giving fundamental definitions and concepts related to the generalized groups, topological generalized groups and sheaves. We can start by giving some basic recalls of the concept of generalized group which was first defined by Molaei

[16].

**Definition 1.** ([16]) A generalized group $G$ is a non-empty set admitting an operation called multiplication subject to the set of rules given below:

(i) $(ab)c = a(bc)$, for all $a, b, c \in G$ (associative law);

(ii) For each $a \in G$, there exists a unique $e(a) \in G$ such that $ae(a) = e(a)a = a$ (existance and uniqueness of identity);

(iii) For each $a \in G$, there exists $a^{-1} \in G$, such that $a\,a^{-1} = a^{-1}\,a = e(a)$ (existance of inverse).

**Example 1.** ([17]) The set $G = \left\{ \begin{bmatrix} a & b \\ c & d \end{bmatrix} : a, b, c \text{ and } d \text{ are real numbers} \right\}$ with the operation

$$\left( \begin{bmatrix} a & b \\ c & d \end{bmatrix}, \begin{bmatrix} e & f \\ g & h \end{bmatrix} \right) \to \begin{bmatrix} a & f \\ g & d \end{bmatrix},$$

is a generalized group in which for all $A \in G$,

$$e(A) = \begin{bmatrix} a & f \\ g & d \end{bmatrix} \text{ and } A^{-1} = \begin{bmatrix} a & f \\ g & d \end{bmatrix},$$

where $e(A)$ and $A^{-1}$ are the identity and the inverse of matrix $A$ respectively.

**Example 2.** ([14]) Let $G = \mathbb{R} \times \{\mathbb{R} \setminus \{0\}\}$. Then with the multiplication $(a, b)(c, d) = (bc, bd)$ is a generalized group in which for all $(a, b) \in G$, $e(a, b) = (a/b, 1)$ and $(a, b)^{-1} = (a/b^2, 1/b)$.

**Example 3.** ([6]) Let $G$ with the multiplication $m$ be generalized group. Then $G \times G$ with the multiplication

$$m_1\big((a, b), (c, d)\big) = \big(m(a, c), m(b, d)\big)$$

is a generalized group. For any element $(a, b) \in G \times G$ the identity element is $e_1(a, b) = \big(e(a), e(b)\big)$ and inverse element is $(a, b)^{-1} = (a^{-1}, b^{-1})$.

**Theorem 1.** ([17]) Each $a$ in a generalized group $G$ has a unique inverse in $G$.

**Example 4.** ([17]) Let $S = \{1,2\}$. Then $S$ with the binary operation: $2.2 = 2, 2.1 = 1.2 = 2, 1.1 = 1$ is a semigroup. Then $S$ is semigroup which is not a generalized group, because the identity of 2 is not unique.

It is easily seen from Definition 1 that every group is a generalized group. But it is not true in general that every generalized group is a group.

**Lemma 1.** ([18]) Let $G$ be a generalized group and $ab = ba$ for all $a, b \in G$. Then, $G$ is an abelian group.

**Definition 2.** ([14]) A non-empty subset $H$ of a generalized group $G$ is a generalized subgroup of $G$ if and only if for all $a, b \in H, ab^{-1} \in H$.

**Theorem 2.** ([17]) Let $G$ be a generalized group such that $a \in G$. Then,

$$G_a = e^{-1}\{e(a)\} = \{x \in G : e(x) = e(a)\}$$

is a generalized subgroup of $G$. Furthermore, $G_a$ is a group.

Let us give some results related to the structure of generalized groups via following lemma.

**Lemma 2.** ([2]) Let $G$ be a generalized group. Then,

(i) $e(a) = e(a^{-1})$ and $e(e(a)) = e(a)$ where $a \in G$.
(ii) $(a^{-1})^{-1} = a$ where $a \in G$.
(iii) The set $\{G_a = e^{-1}\{e(a)\} : a \in G\}$ is a partitation of groups for $G$.

We here state definition of a topological generalized group which was defined by Molaei [16] and set fourth simpliest properties of this structure from topological point of view was presented in [14] and [17] thereof.

**Definition 3.** ([17]) A topological generalized group $G$ is a set which satisfies the following conditions:

(i) $G$ is generalized group;
(ii) $G$ is a Hausdorff topological space;
(iii) The mappings
$$m_1: G \times G \to G, (a,b) \to ab$$
and
$$m_2: G \to G, a \to a^{-1}$$
are continuous mappings.

If $a \in G$ then $G_a = e^{-1}(\{e(a)\})$ with the product of $G$ is a topological group, and $G$ is disjoint union of such topological groups i.e. $G = \vee_{a \in G} G_a$.

**Example 5.** ([17]) Every non-empty Hausdorff topological space $G$ with the operation:
$$m: G \times G \to G$$
$$(a,b) \mapsto a$$
is a topological generalized group.

**Example 6.** ([14]) The set $G = IR \times (IR \setminus \{0\})$ with the topology induced by a Euclidean metric, and with the multiplication $(a,b).(c,d) = (bc,bd)$ is a topological generalized group.

**Definition 4.** ([19]) Let $X, S$ both topological spaces, and $\pi: S \to X$ be a locally topological map. Then the pair $S = (S, \pi)$ or shortly $S$ is called a sheaf over $X$.

In the definition of a sheaf, $X$ is not assumed to satisfy any separation axioms (See in [4]). $S$ is called the sheaf space, $\pi$ the projection map, and $X$ the base space. Let $x$ be an arbitrary point in $X$ and $V$ be an open neighborhood of $x$. A section over $V$ is a continuous map $s: V \to S$ such that $\pi \circ s = id_V$.

Let us denote the collection of all sections of $S$, by $\Gamma(V, S)$ and recall the Whitney sum.

**Definition 5.** ([8,21]) Let $(S_1, \pi_1), (S_2, \pi_2), \ldots, (S_k, \pi_k)$ be sheaves on $X$. Construct product $M_W = \Gamma(W, S_1) \times \Gamma(W, S_2) \times \ldots \times \Gamma(W, S_k)$ for $V, W \subset X$ open sets. Let $\Gamma_V^W \colon M_W \to M_V$ defined by $\Gamma_V^W(s) = (s_1|_V, s_2|_V, \ldots, s_k|_V)$ for $(s_1, s_2, \ldots, s_k) \in M_W$ and $V \subset W$. Then $\{M_W, \Gamma_V^W\}$ is a presheaf. The Whitney sum of $S_1, S_2 \ldots, S_k$ sheaves is a sheaf defined by this presheaf and denoted by $S^* = S_1 \oplus S_2 \oplus \ldots \oplus S_k$.

Now we can say that the Whitney sum of sheaves $(S_1, \pi_1), (S_2, \pi_2), \ldots, (S_k, \pi_k)$:

$$S^* = S_1 \oplus \ldots \oplus S_k := \{\sigma = (\sigma_1, \ldots, \sigma_k) \in S_1 \times \ldots \times S_k : \pi_1(\sigma) = \cdots = \pi_k(\sigma)\}$$

$$= \bigvee_{x \in X} ((S_1)_x \times \ldots \times (S_k)_x)$$

is a set over $X$ topological spaces. Then the map $\pi \colon S^* = S_1 \oplus S_2 \oplus \ldots \oplus S_k \to X$, $\pi(\sigma) = (\pi_i \circ P_i)(\sigma)$ is a local homeomorphism, hence $S^* = S_1 \oplus S_2 \oplus \ldots \oplus S_k$ is a sheaf over $X$.

**Theorem 3.** ([11]) Let $(S_i, \pi_i)$, $i = 1, \ldots, k$ be sheaves and $S^* = S_1 \oplus S_2 \oplus \ldots \oplus S_k$ be Whitney sum of $S_1, S_2 \ldots, S_k$. Then there is bijection $\pi \colon S_x^* \to (S_1)_x \times (S_2)_x \times \ldots \times (S_k)_x$ defined by $(W, (s_1, s_2, \ldots, s_k))_x \to (s_1(x), s_2(x), \ldots, s_k(x))$.

**Theorem 4.** ([11]) Let $(S_i, \pi_i), i = 1, \ldots, k$ be sheaves on $X$. Then the canonic projection $P_i \colon S_1 \oplus S_2 \oplus \ldots \oplus S_k \to S_i$, $P_i(\sigma_1, \sigma_2, \ldots, \sigma_k) = \sigma_i$ is a sheaf morphism.

Let $s_i \in \Gamma(W_i, S_i)$ for $i = 1, \ldots, k$. Define $s_1 \oplus \ldots \oplus s_k \colon W \to S^* = S_1 \oplus S_2 \oplus \ldots \oplus S_k$, such that $(s_1 \oplus \ldots \oplus s_k)(x) = (s_1(x), s_2(x), \ldots, s_k(x))$. Clearly $(s_1, s_2, \ldots, s_k) \in M_W$ and $r(s_1, s_2, \ldots, s_k) = (W, (s_1, s_2, \ldots, s_k))_x = (s_1(x), s_2(x), \ldots, s_k(x)) = (s_1 \oplus \ldots \oplus s_k)(x)$. Therefore since $s_1 \oplus \ldots \oplus s_k = r(s_1, s_2, \ldots, s_k) \in \Gamma(W, S_1 \oplus S_2 \oplus \ldots \oplus S_k)$ we have

$$\Gamma(W, S_1 \oplus S_2 \oplus \ldots \oplus S_k) = \Gamma(W, S_1) \times \Gamma(W, S_2) \times \ldots \times \Gamma(W, S_k).$$

Furthermore since $W \subset X$ is open set and $\pi$ is a local homeomorphism $s(W)$ is open set in $S$, and $S$ is union of these type of open sets. Also if $s_1, s_2 \in \Gamma(W, S)$ and $s_1(x) = s_2(x)$ for $x \in X$ then $s_1 = s_2$ in $W$. So we can say that every element of $S$ can be seen as a substance of sections in S.

## 2. MAIN RESULTS

### 2.1 The Sheaf of The Groups Formed By Topological Generalized Group Over Topological spaces

Let $\mathcal{C}$ be the category of the topological spaces $X$ satisfying the property that all pointed spaces $(X, x)$ with $x \in X$ have same homotopy type. This category includes all topological vector spaces.

Let us take $X \in \mathcal{C}$ as a base set if $P$ is any topological group with identity element $p_0$ as base point. Then the set of homotopy class of homotop maps preserving the base point from $(X, x)$ to $(P, p_0)$ obtained for each $x \in X$, $(X, x)$ pointed topological spaces i. e. $S(X) = \vee_{x \in X}[(X, x), (P, p_0)]$. Thus $S(X)$ is a set over $X$.

If $(P, p_0)$ is any topological group with the identity element of the group is $p_0$, we can construct a sheaf over $X$ by using following theorem which is given by C. Yıldız [23].

**Theorem 5.** ([23]) Let $(P, p_0)$ be any pointed topological group with the identity element $p_0$ and $X \in \mathcal{C}$. If $\pi: S(X) \to X$ such that $\pi(\sigma) = \pi([f]_x) = x$ for $\sigma = [f]_x \in S(X)$, $x \in X$ then there is the natural topology over $S(X)$ such that $\pi$ is locally topological with respect to this topology. Thus the pair $(S, \pi)$ is a sheaf over $X$.

In Theorem 5, C. Yıldız by defining $S(X) = \vee_{x \in X}[(X, x), (P, p_0)]$ and $\pi: S(X) \to X$ such that $\pi(\sigma) = x$, $x \in X$ and a mapping $s: V \to S(X)$ as follows:

If $x_0 \in X$, then there exists a group $[(X, x_0), (P, p_0)]$ in $S(X)$. If $y$ is any point in $V$, Then $(X, x_0)$ and $(X, y)$ are having same homotopy type where $V = V(x_0)$ open neighborhood of $x_0$ in $X$. Therefore, there is a homotopy equivalence map $\Phi: (X, x_0) \to (X, y)$.

Hence from the diagram in Figure 1, the map $h = f \circ \Phi: (X, y) \to (P, p_0)$ is continuous and base-point preserving. $[h]_y \in [(X, y), (P, p_0)]$ is a homotopy class of map $f \circ \Phi = h$.

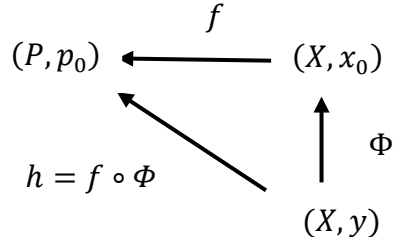

**Figure 1**

Therefore, we define $s(y) = [h]_y$. In this way $s$ is well defined and $(\pi \circ s)(y) = \pi(s(y)) = y$ for each $y \in V$. Therefore $\pi \circ s = I_V$. Thus $s$ is called a section of $S(X)$ over $V$.

Let us denote the collection of all sections of $S(X)$, by $\Gamma(V,S)$. A topology-base is constructed on $S(X)$ by using $s(V) = V_{y \in V}[h]_y$,

$$\beta = \{s(V) : V = V(x) \subset X, x \in X, s \in \Gamma(V,S)\}.$$

Thus gives a natural topology on $S(X)$. Therefore $S(X)$ is a topological space.

Therefore the sheaf $(S, \pi)$ given by Theorem 5 is a sheaf of the homotopic groups formed by topological group $P$ over $(X,x)$ pointed topological spaces [23]. The stalk of the sheaf $(S, \pi)$ over $X$ is the group $[(X,x),(P,p_0)] = \pi^{-1}$ denoted by $S(X)_x$ for every $x \in X$.

$\Gamma(V,S)$ is a group with pointwise multiplication defined by

$$(s_1 s_2)(y) = s_1(y)s_2(y), \; s_1, s_2 \in \Gamma(V,S) \text{ and } y \in V.$$

And in this group the identity element is $I: V \to S$ which is obtained by means of the identity element of $[(X,x),(P,p_0)]$ and the inverse element of $s \in \Gamma(V,S)$ is $s^{-1} \in \Gamma(V,S)$ which is obtained by the inverse element of $[(X,x),(P,p_0)]$. Therefore $(S, \pi)$ is an algebraic sheaf with the operation $(.): S(X) \otimes S(X) \to S(X)$ (that is, $(\sigma_1, \sigma_2) \to \sigma_1 . \sigma_2$ for every $\sigma_1, \sigma_2 \in S(X)$ is continuous [23].

Now let begin to construct the sheaf over $X$ by the finite pointed topological generalized group $P$. We begin with constructing the Whitney sum of sheaves $S_1, \ldots, S_k$  i. e. $S^* = S_1 \oplus \ldots \oplus S_k$.

Let us now define a map $\pi: S^*(X) \to X$, $\pi(\sigma) = (\pi_{i_0} \circ P_{i_0})(\sigma)$ for fixed $i_0$ and canonic projection $P_{i_0}$.

If $x_0 \in X$, then there exists groups $[(X, x_0), (P, p_i)]$ in $S_i(X)$ for $i = 1, \ldots, k$. Let $\sigma = (\sigma_1, \sigma_2, \ldots, \sigma_k) = ([h_1]_{x_0}, [h_2]_{x_0}, \ldots, [h_k]_{x_0})$ be a homotopy class in the group $\prod_{i=1,\ldots,k} [(X, x_0), (P, p_i)]$ where the map $h_i = f_i \circ \Phi: (X, x) \to (P, p_i)$ is continuous base point preserving. $[h_i]_y \in [(X, y), (P, p_i)]$ for $i = 1, \ldots, k$ is a homotopy class of map $f_i \circ \Phi = h_i$.

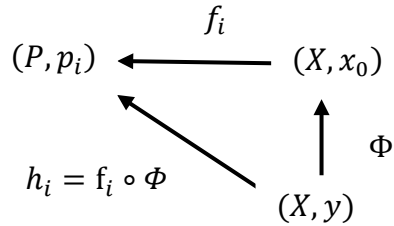

**Figure 2**

If $x_0 \in X$ is an arbitrarily fixed point, then let us denote $V = V(x_0)$ open neighborhood of $x_0$ in $X$. Now, we can define a mapping $s = (s_1, \ldots, s_k): V \to S^*(X)$, as follows:

If $y$ is any point in $V$, then we define $s(y) = (s_1, s_2, \ldots, s_k)(y) = (s_1(y), s_2(y), \ldots, s_k(y))$ for $s_i(y) = [h_i]_y$, $i = 1, \ldots, k$. In this way $s$ is well defined and

1. $(\pi \circ s)(y) = (\pi_i \circ P_i \circ s)(y) = \pi_i \left( P_i(s_1(y), s_2(y), \ldots, s_k(y)) \right) = \pi_i(s_i(y)) = y$ for each $y \in V$. Therefore $\pi \circ s = I_V$.

2. If, $x_0$ is an arbitrary fixed point in $V$,

$$s(x_0) = (s_1, \ldots, s_k)(x_0) = (s_1(x_0), \ldots, s_k(x_0))$$
$$= ([f_1 \circ I_x]_{x_0}, \ldots, [f_k \circ I_x]_{x_0}) = ([f_1]_{x_0}, \ldots, [f_k]_{x_0})$$

for $V = V(x_0)$. Hence it can be written as $s(V)=\prod_{i=1,\ldots,k} s_i(V) = \prod_{i=1,\ldots,k}(\vee_{y \in V} [h_i]_y)$.
If we can define $s(V)$ as an open set, then it can be easily shown that the family

$$\beta = \left\{ s(V) := \prod_{i=1,\ldots,k} s_i(V) : V = V(x) \subset X, x \in X, s_i \in \Gamma(V, S_i) \right\} \cup \{S^*\}$$

is a topology-base on $S^*(X)$. Thus $S^*(X)$ is a topological space.

Now we can show that $\pi: S^*(X) \to X$ is local topological.

If $\sigma = [h]_y \in S^*(X)$ and $y \in X$, then $\pi(\sigma) = \pi([h]_y) = y$. Therefore, there is a map $s: V \to S^*(X)$ such that $s(y) = \sigma$, $y \in V$. Now, let us assume that $U(\sigma) = s(V)$ and $\pi|_U = \pi^*$.

1. The map $\pi^* = \pi|_U : U \to V$ is injective. Because for any $\sigma_1, \sigma_2 \in s(V)$ there are points $y_1, y_2$ respectively in $V$ such that

$$\sigma_1 = s(y_1) = (s_1, s_2, \ldots, s_k)(y_1) = (s_1(y_1), s_2(y_1), \ldots, s_k(y_1))$$
$$= ([f_1 \circ \Phi]_{y_1}, [f_2 \circ \Phi]_{y_1}, \ldots, [\Phi \circ f_k]_{y_1}),$$

$$\sigma_2 = s(y_2) = (s_1, s_2, \ldots, s_k)(y_2) = (s_1(y_2), s_2(y_2), \ldots, s_k(y_2))$$
$$= ([f_1 \circ \Phi']_{y_2}, [f_2 \circ \Phi']_{y_2}, \ldots, [f_k \circ \Phi']_{y_2}).$$

That is, we have the following diagrams for $i = 1, \ldots, k$.

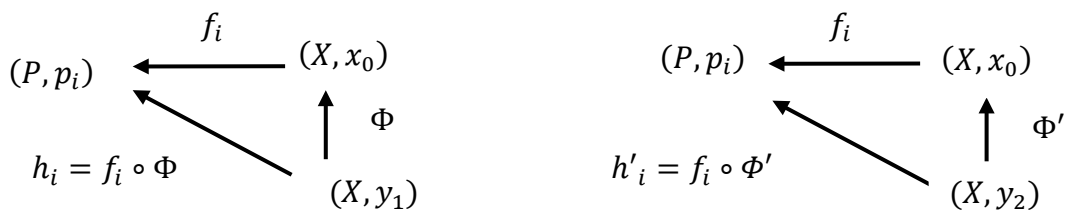

**Figure 3**

If $\pi^*(\sigma_1) = \pi^*(\sigma_2)$, then

$$\pi^*(s(y_1)) = \pi^*(s(y_2)) \Rightarrow \pi^*\left(([f_1 \circ \Phi]_{y_1}, [f_2 \circ \Phi]_{y_1}, \ldots, [f_k \circ \Phi]_{y_1})\right)$$

$$= \pi^*([f_1 \circ \Phi']_{y_2}, [f_2 \circ \Phi']_{y_2}, \ldots, [f_k \circ \Phi']_{y_2}) \Rightarrow y_1 = y_2.$$

Since for each $i = 1, \ldots, k$,

$$\Phi \sim \Phi' \Rightarrow f_i \circ \Phi \sim f_i \circ \Phi' \Rightarrow [f_i \circ \Phi]_{y_1} = [f_i \circ \Phi']_{y_2} \Rightarrow \sigma_1 = \sigma_2.$$

2. The map $\pi^* = \pi|_U : U \to V$ is continuous. In fact, if $\sigma \in U = s(V) \Rightarrow \pi^*(\sigma) = y \in V$ and $W = W_y \subset V$ is neighbourhood of $y$, then $s(W) \subset U = s(V)$ is neighborhood of $\sigma$ and $\pi^*(s(W)) = W \subset V$. So $\pi^*$ is continuous.

3. $\pi^{*-1} = (\pi|_U)^{-1} = s : V \to U = s(V)$ is continuous. In fact, if $y$ is any point in $V$, $s(y) = \sigma \in U$ and $U' = U'(\sigma) \subset U$ is a neighborhood of $y$ in $V$ and $s(\pi|_U)(U') \subset U$. So $\pi^{*-1}$ is continuous.

Therefore $\pi$ is locally topological map. Now we can give the following theorem can be given.

**Theorem 6.** Let $(P, p_i)_{i=1,\ldots,k}$ be any pointed finite topological generalized group with the identity elements $p_1, p_2, \ldots, p_k$ and $X \in \mathcal{C}$. If

$$S^* = S_1 \oplus S_2 \oplus \ldots \oplus S_k \text{ and } \pi : S^*(X) \to X$$

such that

$$\pi(\sigma) = (\pi_i \circ P_i)([h_1]_x, [h_2]_x, \ldots, [h_k]_x) = x, \ i = 1, \ldots, k,$$

for $\sigma \in S^*(X)$ and $x \in X$, then there is the natural topology based on $S^*(X)$, such that $\pi$ is locally topological with respect to this natural topology. Thus the pair $(S^*, \pi)$ is a sheaf over $X$.

**Definition 6.** The sheaf $(S^*, \pi)$ given by Theorem 6 is called sheaf of the homotopic groups formed by generalized groups over $(X, x)$, $x \in X$ pointed topological spaces.

**Definition 7.** The group $\prod_{i=1,\ldots,k}[(X,x),(P,p_i)] = \pi^{-1}$ is called the stalk of the sheaf $(S^*, \pi)$ over $X$ and denoted by $S^*(X)_x$ for every $x \in X$.

Now, if $x \in X$ is an arbitrarily fixed point and $V$ is open neighborhood of $x$ in $X$, the mapping $s: V \to S^*(X)$ as defined in the construction of topology of $S^*(X)$, is called section of $S^*(X)$, over $V$. Let us denote the collection of all sections of $S^*(X)$, by $\Gamma(V, S^*)$.

**Theorem 7.** $\Gamma(V, S^*)$ is a group with the operation

$$(s_1 s_2)(y) = s_1(y) s_2(y), s_1, s_2 \in \Gamma(V, S^*)$$

where $y \in V$.

**Proof.** If we consider pointwise multiplication $(s_1^i, s_2^i)(y) = s_1^i(y) s_2^i(y), s_1^i, s_2^i \in \Gamma(V_i, S_i)$ and $y \in V_i$ which is defined on $\in \Gamma(V_i, S_i)$ for $i = 1, \ldots, k$. Proof follows from that the operation of production is well-defined and closed. Clearly, the operation of production is associative and the mapping $I: V \to S^*(X)$ is identity element which is obtained by means of the identity element of $\prod_{i \in I}[(X,x),(P,p_i)]$. On the other hand, the any inverse element of $s \in \Gamma(V, S^*)$, namely, $s^{-1} \in \Gamma(V, S^*)$ which is obtained by means of the homotopy inverses of pointed groups $(P, p_i)$ for $i = 1, \ldots, k$. Hence $\Gamma(V, S^*)$ is a group.

From the Theorem 6, $(S^*, \pi)$ is an algebraic sheaf with the continuous operation

$$(.): S^*(X) \otimes S^*(X) \to S^*(X),$$
$$(\sigma_1, \sigma_2) \to \sigma_1 . \sigma_2$$

where $\sigma_1, \sigma_2 \in S^*(X)$.

**CONFLICT OF INTEREST**

No conflict of interest was declared by the authors.